\documentclass[reqno,12pt]{amsart}
\usepackage{amsfonts}
\usepackage{amsmath}
\usepackage{graphicx}
\usepackage{amscd}

\newtheorem{theorem}{Theorem}
\theoremstyle{plain}
\newtheorem{acknowledgement}{Acknowledgement}

\newtheorem{corollary}{Corollary}

\numberwithin{equation}{section}

\input{tcilatex}

\begin{document}
\author{}
\title{}
\maketitle


\begin{center}
\thispagestyle{empty}\textbf{{\Large {A NOTE ON p-ADIC q-EULER MEASURE}}}

\bigskip 

\textbf{{\Large {Hacer OZDEN}}}$^{\ast }$\textbf{{\Large {, Yilmaz\ SIMSEK}}}%
$^{^{\ast \ast }}$\textbf{{\Large {, Seog-Hoon RIM}}}$^{\ast \ast \ast }$%
\textbf{{\Large {, Ismail Naci CANGUL}}}$^{\ast }$

\bigskip 

$^{\ast }${\Large \textit{{University of Uludag, F}aculty of Arts and
Science, Department of Mathematics, 16059 Bursa}, Turkey}

$^{\ast \ast }$\textit{{\Large {University of Akdeniz, Faculty of Arts and
Science, Department of Mathematics, }}}{\Large {\textit{07058 Antalya},
Turkey}}

\bigskip 

$^{\ast \ast \ast }$\textit{{\Large {University of Kyungpook, Department of
Mathematical Education, Taegu, 702-701, S. KOREA}}}

\bigskip 

{\Large 
}

{\Large \textbf{E-Mail: hozden@uludag.edu.tr, simsek@akdeniz.edu.tr,
shrim@knu.ac.kr and cangul@uludag.edu.tr}}

{\Large \ }

\textbf{{\Large {Abstract}}}
\end{center}

In this paper, we will investigate some interesting properties of the
modified q-Euler numbers and polynomials. The main purpose of this paper is
to construct p-adic q-Euler measure on $%
\mathbb{Z}
_{p}$.



\bigskip

\section{Introduction, Definitions and Notations}

Let $p$ be a fixed odd prime. Throughout this paper $%
\mathbb{Z}
_{p},~%
\mathbb{Q}
_{p},~%
\mathbb{C}
~$and $%
\mathbb{C}
_{p}$ will respectively denote the ring of $p$-adic rational integers, the
field of $p$-adic rational numbers, the complex number field and the
completion of the algebraic closure of $%
\mathbb{Q}
_{p}$. Let $v_{p}$ be the normalized exponential valuation of $%
\mathbb{C}
_{p}$ with $|p|_{p}~=~p^{-v_{p}(p)}~=~1/p,~$see~\cite{kimmodified}, \cite%
{kimon a q},~\cite{kimon the q}. When one talks of $q$-extensions, $q~$is
variously considered as an indeterminate, a complex $q$ $\in 
\mathbb{C}
$,~or a $p$-adic number $q$ $\in 
\mathbb{C}
_{p}.$ If $q$ $\in 
\mathbb{C}
$,~one normally assumes $|q|~<~1$. If $q$ $\in 
\mathbb{C}
_{p},$ then we assume $|q-1|~<~p^{-\frac{1}{p-1}}$, so that $q^{x}=~\exp
(x\log q)$ for $|x|_{p}~\leq $ 1, see \cite{kimsums}, \cite{kiman inv},~\cite%
{srikimsim}. It was well known that Euler numbers are defined by

\begin{equation*}
\frac{2}{e^{t}+1}=\overset{\infty }{\underset{n=0}{\dsum }}E_{n}\frac{t^{n}}{%
n!},\text{~see \cite{kimmodified},~\cite{kimon a q},~\cite{kimon the q},~%
\cite{kimsums},~\cite{kiman inv},~\cite{srikimsim}.}
\end{equation*}

In the recent paper \cite{kimon the q}, the q-extension of Euler numbers are
defined inductively by

\begin{equation*}
E_{0,q}=1,~q(qE+1)^{n}+E_{n,q}=\left\{ 
\begin{array}{cc}
\lbrack 2]_{q} & if~n=0 \\ 
0 & if~n\neq 0%
\end{array}%
\right.
\end{equation*}%
with the usual convention of replacing $E^{m}~$by $E_{m,q}.$ In \cite%
{kimmodified}, the definition of modified $q$-Euler numbers $\mathcal{E}%
_{0,q}~$is introduced by

\begin{equation*}
\mathcal{E}_{0,q}=\frac{[2]_{q}}{2},~(q\mathcal{E}+1)^{k}-\mathcal{E}%
_{k,q}=\left\{ 
\begin{array}{cc}
\lbrack 2]_{q} & if~k=0 \\ 
0 & if~k>0%
\end{array}%
\right.
\end{equation*}%
with the usual convention of replacing $\mathcal{E}^{i}~$by $\mathcal{E}%
_{i,q}.$

For a fixed positive integer $d$ with $(p,d)=1$, set

\begin{equation*}
X_{d}=\lim_{\overset{\leftarrow }{N}}%
\mathbb{Z}
/dp^{N}%
\mathbb{Z}
,\ \ X_{1}=%
\mathbb{Z}
_{p,}
\end{equation*}

\begin{equation*}
X^{\ast }=\underset{\underset{(a,p)=1}{0<a<dp}}{\cup ~}(a+dp%
\mathbb{Z}
_{p}),
\end{equation*}

\begin{equation*}
a+dp^{N}%
\mathbb{Z}
_{p}=\left\{ x\in ~X~:~x\equiv a\ (\func{mod}~dp^{N})\right\} ,
\end{equation*}%
where $a\in $\bigskip\ $%
\mathbb{Z}
$ satisfies the condition $0\leq a<dp^{N},$ see \cite{srikimsim}.

We say that $f$ is a uniformly differentiable function at a point $a\in $%
\bigskip\ $%
\mathbb{Z}
_{p}$, and write $f\in $\bigskip\ $UD(%
\mathbb{Z}
_{p}),$ if the difference quotient

\begin{equation*}
F_{f}(x,y)=\frac{f(x)-f(y)}{x-y}
\end{equation*}%
has a limit $f%
{\acute{}}%
(a)$ as $(x,y)\rightarrow (a,a).$ For~ $f\in $\bigskip\ $UD(%
\mathbb{Z}
_{p}),$ an invariant $p$-adic $q$-integral was defined by

\begin{equation*}
I_{q}(f)=\dint\limits_{%
\mathbb{Z}
_{p}}f(x)d\mu _{q}(x)=\lim_{N\rightarrow \infty }\frac{1}{[p^{N}]_{q}}%
\overset{p^{N}-1}{\underset{x=0}{\dsum }}f(x)g^{x},~\text{see \cite{kim anew}%
.}
\end{equation*}

\bigskip The $q$-extension of $n\in 
\mathbb{N}
$ is defined by

\begin{equation*}
\lbrack n]_{q}=\frac{1-q^{n}}{1-q}=1+q+q^{2}+...+q^{n-1},
\end{equation*}%
and

\begin{equation*}
\lbrack n]_{-q}=\frac{1-(-q)^{n}}{1-(-q)}=1-q+q^{2}-...+(-q)^{n-1},~\text{%
see \cite{kimon the q},\cite{kimsums},\cite{kiman inv}.}
\end{equation*}

The modified $p$-adic $q$-integral on $%
\mathbb{Z}
_{p}$ is defined by

\begin{equation*}
I_{-q}(f)=\dint\limits_{%
\mathbb{Z}
_{p}}f(x)d\mu _{-q}(x),
\end{equation*}%
where $d\mu _{-q}(x)$ $=\underset{q\rightarrow -q}{\lim }d\mu _{q}(x)$. In
this paper, we will investigate some interesting properties of the modified $%
q$-Euler numbers and polynomials. The purpose of this paper is to construct $%
p$-adic $q$-Euler measure on $%
\mathbb{Z}
_{p}$.

\section{\protect\bigskip $\mathbf{p}$\textbf{-adic }$\mathbf{q}$\textbf{%
-Euler Measure on }$%
\mathbb{Z}
_{p}$}

$q$-Euler numbers are known by

\begin{equation*}
\mathcal{E}_{0,q}=\frac{[2]_{q}}{2},~(q\mathcal{E}+1)^{n}-\mathcal{E}%
_{n,q}=\left\{ 
\begin{array}{cc}
\lbrack 2]_{q} & if~n=0 \\ 
0 & if~n>0%
\end{array}%
\right.
\end{equation*}%
with the usual convention of replacing $\mathcal{E}^{i}~$by $\mathcal{E}%
_{i,q},~$see \cite{kimmodified}. It was known that the $q$-Euler numbers can
be represented by $p$-adic $q$-integrals on $%
\mathbb{Z}
_{p}$ as follows:

\begin{equation*}
\mathcal{E}_{n,q}=\dint\limits_{%
\mathbb{Z}
_{p}}q^{-x}[x]_{q}^{n}d\mu _{-q}(x)=[2]_{q}\left( \frac{1}{1-q}\right) ^{n}%
\overset{n}{\underset{l=0}{\dsum }}(-1)^{l}\left( 
\begin{array}{c}
n \\ 
l%
\end{array}%
\right) \frac{1}{1+q^{l}},~\text{see \cite{kimmodified}}.
\end{equation*}

We now also consider the $q$-Euler polynomials as follows:

\begin{equation*}
\mathcal{E}_{n,q}(x)=\dint\limits_{%
\mathbb{Z}
_{p}}q^{-x}[t+x]_{q}^{n}d\mu _{-q}(t)=~\overset{n}{\underset{l=0}{\dsum }}%
\left( 
\begin{array}{c}
n \\ 
l%
\end{array}%
\right) q^{xl}\mathcal{E}_{l,q}[x]_{q}^{n-l}.
\end{equation*}%
Thus, we note that

\begin{equation}
\mathcal{E}_{n,q}(x)=[d]_{q}^{n}\frac{[2]_{q}}{[2]_{q^{d}}}\overset{d-1}{%
\underset{a=0}{\dsum }}\left( -1\right) ^{a}\mathcal{E}_{n,q^{d}}\left( 
\frac{n+a}{d}\right) ,~\text{see \cite{kimmodified}.}  \label{1}
\end{equation}

Let $\chi $ be the Dirichlet's character with odd conductor $d\in 
\mathbb{N}
$, and let $F_{\chi ,q}(t)$ be the generating function of $\mathcal{E}%
_{n,\chi ,q}$ as follows:

\begin{equation}
F_{\chi ,q}(t)=[2]_{q}\overset{\infty }{\underset{n=0}{\dsum }}\left(
-1\right) ^{n}\chi (n)e^{\left[ n\right] }q^{t}=\overset{\infty }{\underset{%
n=0}{\dsum }}\mathcal{E}_{n,\chi ,q}\frac{t^{n}}{n!}\text{,~\bigskip see 
\cite{kimmodified}.}  \label{2}
\end{equation}%
From (\ref{2}) we derive

\begin{equation}
\mathcal{E}_{n,\chi ,q}=[d]_{q}^{n}\frac{[2]_{q}}{[2]_{q^{d}}}\overset{d-1}{%
\underset{a=0}{\dsum }}\left( -1\right) ^{a}\chi (a)\mathcal{E}%
_{n,q^{d}}\left( \frac{a}{d}\right) .  \label{3}
\end{equation}

For any positive integers $N,~k$ and $d~$(odd), let $\mu _{k}^{\ast }=\mu
_{k,q;E}^{\ast }$ be defined by

\begin{equation}
\mu _{k}^{\ast }(a+dp^{N}%
\mathbb{Z}
_{p})=\left( -1\right) ^{a}[dp^{N}]_{q}^{k}\frac{[2]_{q}}{[2]_{q^{dp^{N}}}}%
\mathcal{E}_{k,q^{dp^{N}}}\left( \frac{a}{dp^{N}}\right) .  \label{4}
\end{equation}

Then we see that

\begin{eqnarray*}
&&\overset{p-1}{\underset{i=0}{\dsum }}\mu _{k}^{\ast }(a+idp^{N}+dp^{N+1}%
\mathbb{Z}
_{p}) \\
&=&[dp^{N+1}]_{q}^{k}\frac{[2]_{q}}{[2]_{q^{dp^{N+1}}}}\overset{p-1}{%
\underset{i=0}{\dsum }}\left( -1\right) ^{a+idp^{N}}\mathcal{E}%
_{k,q^{dp^{N+1}}}\left( \frac{a+idp^{N}}{dp^{N+1}}\right)  \\
&=&\left( -1\right) ^{a}[dp^{N+1}]_{q}^{k}\frac{[2]_{q}}{[2]_{q^{dp^{N+1}}}}%
\overset{p-1}{\underset{i=0}{\dsum }}\left( -1\right) ^{i}\mathcal{E}%
_{k,(q^{dp^{N}})^{p}}\left( \frac{\frac{a}{dp^{N}}+i}{p}\right)  \\
&=&\left( -1\right) ^{a}[dp^{N}]_{q}^{k}\frac{[2]_{q}}{[2]_{q^{dp^{N}}}}%
\left( \frac{[2]_{q^{dp^{N}}}}{[2]_{(q^{dp^{N}})^{p}}}[p]_{q^{dp^{N}}}^{k}%
\overset{p-1}{\underset{i=0}{\dsum }}\left( -1\right) ^{a}\mathcal{E}%
_{k,(q^{dp^{N}})^{p}}\left( \frac{\frac{a}{dp^{N}}+i}{p}\right) \right)  \\
&=&\left( -1\right) ^{a}[dp^{N}]_{q}^{k}\frac{[2]_{q}}{[2]_{q^{dp^{N}}}}%
\mathcal{E}_{k,q^{dp^{N}}}\left( \frac{a}{dp^{N}}\right)  \\
&=&\mu _{k}^{\ast }(a+dp^{N}%
\mathbb{Z}
_{p}).
\end{eqnarray*}%
\bigskip It is easy to see that $|\mu _{k}^{\ast }|~\leq ~M,$ for some
constant $M$. Therefore we obtain the following:

\begin{theorem}
\label{Th 1}For any positive integers $N,~k$ and $d~$(odd), let $\mu
_{k}^{\ast }=\mu _{k,q;E}^{\ast }$ be defined by%
\begin{equation*}
\mu _{k}^{\ast }(a+dp^{N}%
\mathbb{Z}
_{p})=.\left( -1\right) ^{a}[dp^{N}]_{q}^{k}\frac{[2]_{q}}{[2]_{q^{dp^{N}}}}%
\varepsilon _{k,q^{dp^{N}}}\left( \frac{a}{dp^{N}}\right) .
\end{equation*}%
Then $\mu _{k}^{\ast }$ is a measure on $X$.
\end{theorem}

\bigskip From the definition of $\mu _{k}^{\ast }$, we derive the following:

\begin{eqnarray}
&&\dint\limits_{X}\chi (x)d\mu _{_{k}}^{\ast }(x)=\lim_{N\rightarrow \infty
}~\overset{dp^{N}-1}{\underset{a=0}{\dsum }}\chi (a)\mu _{_{k}}^{\ast
}(a+dp^{N}%
\mathbb{Z}
_{p})  \notag \\
&&=\lim_{N\rightarrow \infty }\overset{dp^{N}-1}{\underset{x=0}{\dsum }}%
(-1)^{x}\chi (x)[dp^{N}]_{q}^{k}\frac{[2]_{q}}{[2]_{q^{dp^{N}}}}\mathcal{E}%
_{k,q^{dp^{N}}}\left( \frac{x}{dp^{N}}\right)   \notag \\
&&=\lim_{N\rightarrow \infty }[d]_{q}^{k}[p^{N}]_{q^{d}}^{k}\frac{[2]_{q}}{%
[2]_{q^{d}}}\frac{[2]_{q^{d}}}{[2]_{(q^{d})^{p^{N}}}}\overset{dp^{N}-1}{%
\underset{x=0}{\dsum }}(-1)^{x}\chi (x)\mathcal{E}_{k,q^{dp^{N}}}\left( 
\frac{x}{dp^{N}}\right)   \label{1a} \\
&=&[d]_{q}^{k}\frac{[2]_{q}}{[2]_{q^{d}}}\lim_{N\rightarrow \infty }\frac{%
[2]_{q^{d}}}{[2]_{(q^{d})^{p^{N}}}}[p^{N}]_{q^{d}}^{k}\overset{d-1}{\underset%
{a=0}{\dsum }}\overset{p^{N}-1}{\underset{x=0}{\dsum }}(-1)^{a+dx}\chi (a+dx)%
\mathcal{E}_{k,(q^{d})^{p^{N}}}\left( \frac{\frac{a}{d}+x}{p^{N}}\right)  
\notag \\
&=&[d]_{q}^{k}\frac{[2]_{q}}{[2]_{q^{d}}}\lim_{N\rightarrow \infty }\left\{ 
\overset{d-1}{\underset{a=0}{\dsum }}\frac{[2]_{q^{d}}}{[2]_{(q^{d})^{p^{N}}}%
}[p^{N}]_{q^{d}}^{k}\overset{p^{N}-1}{\underset{x=0}{\dsum }}(-1)^{a}\chi (a)%
\mathcal{E}_{k,(q^{d})^{p^{N}}}\left( \frac{\frac{a}{d}+x}{p^{N}}\right)
\right\}   \notag
\end{eqnarray}%
and by equation (\ref{1}), this becomes

\begin{equation*}
\dint\limits_{X}\chi (x)d\mu _{_{k}}^{\ast }(x)=[d]_{q}^{k}\frac{[2]_{q}}{%
[2]_{q^{d}}}\overset{d-1}{\underset{a=0}{\dsum }}(-1)^{a}\chi (a)\mathcal{E}%
_{n,q^{d}}\left( \frac{a}{d}\right) .
\end{equation*}%
By using (\ref{3}) and (\ref{1a}), we obtain the following:

\begin{theorem}
\label{Th 2}For any positive integer $k$, we have%
\begin{equation*}
\dint\limits_{X}\chi (x)d\mu _{_{k}}^{\ast }(x)=\mathcal{E}_{k,\chi ,q.}
\end{equation*}
\end{theorem}

From Theorem \ref{Th 1} and equation (\ref{1}), we note that

\begin{eqnarray*}
&&\mu _{k}^{\ast }(a+dp^{N}%
\mathbb{Z}
_{p}) \\
&=&\left( -1\right) ^{a}\frac{[2]_{q}}{[2]_{q^{dp^{N}}}}[dp^{N}]_{q}^{k}%
\mathcal{E}_{k,q^{dp^{N}}}\left( \frac{a}{dp^{N}}\right)  \\
&=&\left( \overset{k}{\underset{l=0}{\dsum }}\left( 
\begin{array}{c}
k \\ 
l%
\end{array}%
\right) q^{dp^{N}\frac{a}{dp^{N}}l}\mathcal{E}_{l,q^{dp^{N}}}[\frac{a}{dp^{N}%
}]_{qdp^{N}}^{k-l}\right) (-1)^{a}\frac{[2]_{q}}{[2]_{q^{dp^{N}}}}%
[dp^{N}]_{q}^{k} \\
&=&\mathcal{E}_{0,q}[\frac{a}{dp^{N}}]_{qdp^{N}}^{k}(-1)^{a}\frac{[2]_{q}}{%
[2]_{q^{dp^{N}}}}[dp^{N}]_{q}^{k}+(-1)^{a}\frac{[2]_{q}}{[2]_{q^{dp^{N}}}}%
[dp^{N}]_{q}^{k}\overset{}{\underset{}{\overset{k}{\underset{l=1}{\dsum }}%
\left( 
\begin{array}{c}
k \\ 
l%
\end{array}%
\right) q^{l}}} \\
&=&\mathcal{E}_{0,q^{dp^{N}}}[\frac{a}{dp^{N}}]_{qdp^{N}}^{k}(-1)^{a}\frac{%
[2]_{q}}{[2]_{q^{dp^{N}}}}[dp^{N}]_{q}^{k}+(-1)^{a}\frac{[2]_{q}}{%
[2]_{q^{dp^{N}}}}[dp^{N}]_{q}^{k}\cdot  \\
&&\overset{}{\underset{}{\overset{k}{\underset{l=1}{\dsum }}\left( 
\begin{array}{c}
k \\ 
l%
\end{array}%
\right) q^{l}\mathcal{E}_{l,q^{dp^{N}}}\left[ \frac{a}{dp^{N}}\right]
_{q^{dp^{N}}}^{k-l}}} \\
&=&\frac{[2]_{q^{dp^{n}}}}{2}\frac{[a]_{q}^{k}}{[dp^{N}]_{q}^{k}}(-1)^{a}%
\frac{[2]_{q}}{[2]_{q^{dp^{N}}}}[dp^{N}]_{q}^{k}+(-1)^{a}\frac{[2]_{q}}{%
[2]_{q^{dp^{N}}}}[dp^{N}]_{q}^{k}\cdot . \\
&&~\overset{k}{\underset{l=1}{\dsum }}\left( 
\begin{array}{c}
k \\ 
l%
\end{array}%
\right) q^{l}\mathcal{E}_{l,q^{dp^{N}}}\left( \frac{\left[ a\right] _{q}}{%
\left[ dp^{N}\right] _{q}}\right) ^{k-l}.
\end{eqnarray*}

Thus, we have

\begin{eqnarray*}
\lim_{N\rightarrow \infty }~\mu _{k}^{\ast }(a+dp^{N}%
\mathbb{Z}
_{p}) &=&\frac{[2]_{q}}{2}(-1)^{a}[a]_{q}^{k} \\
&=&q^{a}\frac{[2]_{q}}{2}(-1)^{a}q^{-a}[a]_{q}^{k} \\
&=&q^{-a}[a]_{q}^{k}\lim_{N\rightarrow \infty }~\mu _{-q}(a+dp^{N}%
\mathbb{Z}
_{p}),
\end{eqnarray*}%
where $d$ is a positive odd integer. Therefore we obtain the following:

\begin{theorem}
\label{Th 3}For any positive integer $k$, we have%
\begin{equation*}
q^{-x}[x]_{q}^{k}d\mu _{-q}^{{}}(x)=d\mu _{k}^{\ast }(x).
\end{equation*}
\end{theorem}

\begin{corollary}
\label{Co 1}Let $k$ be a positive integer. Then we have%
\begin{equation*}
\mathcal{E}_{k,\chi ,q}=\dint\limits_{X}\chi (x)d\mu _{_{k}}^{\ast
}(x)=\dint\limits_{X}\chi (x)q^{-x}[x]_{q}^{k}d\mu _{_{-q}}(x).
\end{equation*}
\end{corollary}

Let $\overline{d}$\bigskip\ $=(d,p)~$be the least common multiple of the
conductor $d$ of $X~$and $p$, and let $\mathcal{E}_{n,\chi ,q}$ denote the $n
$-th generalized $q$-Euler number belonging to the character $\chi .~$Then \
we have the $\NEG{q}$-analogue form of Witt's formula in the cyclotomic
field $%
\mathbb{Q}
_{p}(\chi )$ as follows:

For all $n\geq 0,$ we have

\begin{equation}
\mathcal{E}_{n,\chi ,q}=\lim_{\rho \rightarrow \infty }~~\frac{[2]_{q}}{2}%
\overset{\overline{d}\bigskip p^{\rho }\ }{\underset{x=1}{\dsum }}\left(
-1\right) ^{\chi }\chi (x)[x]_{q}^{n}.  \label{6}
\end{equation}

Herein as usual we set $\chi (x)=0$ if $x$ is not prime to the conductor $d$%
. From (\ref{6}) we derive

\begin{eqnarray*}
\mathcal{E}_{n,\chi ,q} &=&\lim_{\rho \rightarrow \infty }~\frac{[2]_{q}}{2}%
\overset{}{\underset{1\leq x\leq \overline{d}\bigskip p^{\rho }\ }{%
\dsum^{\ast }}}\left( -1\right) ^{x}\chi (x)[x]_{q}^{n} \\
&&~~~\ \ \ \ \ \ \ \ \ \ \ +\lim_{\rho \rightarrow \infty }~~\frac{[2]_{q}}{%
\left[ 2\right] _{q^{p}}}~\frac{[2]_{q^{p}}}{2}\overset{}{\underset{1\leq
x\leq \overline{d}\bigskip p^{\rho }\ }{\dsum^{\ast }}}\left( -1\right)
^{x}\chi (p)\chi (x)[p]_{q}^{n}[x]_{q^{p}}^{n} \\
&=&\lim_{\rho \rightarrow \infty }~~\frac{[2]_{q}}{2}\overset{}{\underset{%
1\leq x\leq \overline{d}\bigskip p^{\rho }\ }{\dsum^{\ast }}}\left(
-1\right) ^{x}\chi (x)[x]_{q}^{n} \\
&&~~\ \ \ \ \ \ \ \ \ \ \ \ +\frac{[2]_{q}}{\left[ 2\right] _{q^{p}}}%
~\lim_{\rho \rightarrow \infty }~~\frac{[2]_{q^{p}}}{2}\overset{}{\chi
(p)[p]_{q}^{n}\underset{1\leq x\leq \overline{d}\bigskip p^{\rho }\ }{%
\dsum^{\ast }}}\left( -1\right) ^{x}\chi (x)[x]_{q^{p}}^{n},
\end{eqnarray*}%
where $\ast $means taking the sum over the rational integers prime to $p~$in
the given range. Thus we have

\begin{equation*}
\mathcal{E}_{n,\chi ,q}=~\lim_{\rho \rightarrow \infty }~~\frac{[2]_{q}}{2}%
\overset{}{\underset{1\leq x\leq \overline{d}\bigskip p^{\rho }\ }{%
\dsum^{\ast }}}\left( -1\right) ^{\chi }\chi (x)[x]_{q}^{n}+\frac{[2]_{q}}{%
\left[ 2\right] _{q^{p}}}\chi (p)[p]_{q}^{n}~\mathcal{E}_{n,\chi ,q^{p},}
\end{equation*}%
that is,

\begin{eqnarray*}
\mathcal{E}_{n,\chi ,q}-~\frac{[2]_{q}}{\left[ 2\right] _{q^{p}}}\chi
(p)[p]_{q}^{n}~\mathcal{E}_{n,\chi ,q^{p}} &=&\lim_{\rho \rightarrow \infty
}~~\frac{[2]_{q}}{2}\overset{}{\underset{1\leq x\leq \overline{d}\bigskip
p^{\rho }\ }{\dsum^{\ast }}}\left( -1\right) ^{\chi }\chi (x)[x]_{q}^{n} \\
&=&\dint\limits_{X^{\ast }}(-1)^{x}\chi w^{n}(x)\left\langle x\right\rangle
_{q}^{n}~d\mu _{_{-q}}(x),
\end{eqnarray*}%
where $\left\langle x\right\rangle _{q}=\frac{[x]_{q}}{w(x)},$ and $w(x)$ is
the Teichm\"{u}ller character.

\begin{acknowledgement}
The first and fourth authors are supported by the research fund of Uludag
University project no: F-2006/40. The second author is supported by the
research fund of Akdeniz University. 
\end{acknowledgement}

\end{document}